\documentclass[reqno,11pt]{amsart}
\usepackage[english]{babel}
\usepackage[utf8]{inputenc}
\usepackage[T1]{fontenc}
\usepackage[left=2cm,right=2cm,top=3.5cm,bottom=3cm]{geometry}
\usepackage[shortlabels]{enumitem}
\usepackage{listings}
\usepackage{csquotes}
\usepackage{mathtools}
\usepackage{hyperref}
\usepackage{cleveref}
\usepackage{xcolor}
\usepackage{tikz}

\usetikzlibrary{decorations.markings}
\tikzset{%
	myvertex/.style = {circle, draw, fill, very thick, outer sep=3.5pt, inner sep=0pt, minimum size=#1},
	myvertex/.default = 3.5pt,
	myblob/.style = {myvertex, fill=white, minimum size=#1},
	myedge/.style = {draw=black, very thick, line cap=round, line join=round},
	mysquare/.style = {rectangle, draw, fill, thick, outer sep=3.5pt, inner sep=3.5pt, fill=white, minimum size=#1},
}

\definecolor{ugentblue}{RGB}{30,100,200}
\definecolor{ugentyellow}{RGB}{120,190,0}
\definecolor{ugentred}{RGB}{220,78,40}

\title{Symbolic computation in cubic Jordan matrix algebras and in related structures}
\author{Torben Wiedemann}
\address{RPTU University Kaiserslautern-Landau, Department of Mathematics, Gottlieb-Daimler-Straße 48, 67663 Kaiserslautern, Germany}


\makeatletter
\g@addto@macro\bfseries{\boldmath}
\makeatother

\newlength{\fixmidfigure}
\newcommand{\premidfigure}{\setlength{\fixmidfigure}{\lastskip}\addvspace{-\lastskip}}
\newcommand{\postmidfigure}{\addvspace{\fixmidfigure}}

\newcommand{\cjma}{\texttt{CubicJordanMatrixAlg}}

\newtheoremstyle{default-lem}{}{}{\itshape}{}{\bfseries}{.}{ }{\thmnumber{#2} \thmname{#1}\thmnote{ (#3)}}
\theoremstyle{default-lem}
\newtheorem{lemma}{Lemma}[section]

\newtheoremstyle{default-lem}{}{}{}{}{\bfseries}{.}{ }{\thmnumber{#2} \thmname{#1}\thmnote{ (#3)}}
\theoremstyle{default-lem}
\newtheorem{definition}[lemma]{Definition}
\newtheorem{notation}[lemma]{Notation}
\newtheorem{convention}[lemma]{Convention}
\newtheorem{example}[lemma]{Example}
\newtheorem{remark}[lemma]{Remark}

\crefname{remark}{Remark}{Remarks}
\crefname{example}{Example}{Examples}
\crefname{convention}{Convention}{Conventions}
\crefname{definition}{Definition}{Definitions}
\crefname{notation}{Notation}{Notations}
\crefname{section}{Section}{Sections}
\crefname{figure}{Figure}{Figures}

\newtheoremstyle{noname}{}{}{}{}{\bfseries}{.}{ }{\thmnumber{#2}\thmnote{ #3}}
\theoremstyle{noname}
\newtheorem{miscthm}[lemma]{}

\newcommand{\defl}{\coloneq}
\newcommand{\maparg}[1]{\ifblank{#1}{\:\cdot\:}{#1}}
\newcommand*{\map}[5]{\ifblank{#1}{}{#1\colon} \ifblank{#2}{}{#2 \rightarrow #3} \ifblank{#4}{}{\ifblank{#2}{}{,} #4 \mapsto #5}}
\DeclareMathOperator{\End}{End}
\newcommand{\IN}{\mathbb{N}}
\DeclareMathOperator{\Aut}{Aut}

\newcommand{\roots}{\Phi}
\newcommand{\reflbr}[1]{\sigma(#1)}

\newcommand{\Fgrad}[1]{(F_{4})_{#1}}
\newcommand{\Ggrad}[1]{(G_{2})_{#1}}
\newcommand{\rootbase}{\Delta}

\newcommand{\alt}{A}
\newcommand{\conic}{C}
\newcommand{\comring}{k}
\newcommand{\connorm}{n}
\DeclareMathOperator{\contr}{tr}
\newcommand{\conconj}[1]{\ifblank{#1}{\overline{\:\cdot\:}}{\overline{#1}}}

\DeclareMathOperator{\Her}{Her}
\newcommand{\jor}{J}
\newcommand{\perm}{\mathcal{P}}
\newcommand{\cubel}[2]{#1[#2]}
\newcommand{\cubnorm}{N}
\newcommand{\cubtimes}{\times}
\newcommand{\cubadj}{\sharp}
\newcommand{\cubtr}{T}
\newcommand{\tribrace}[3]{\{\: \maparg{#1},\maparg{#2},\maparg{#3}\:\}}

\newcommand{\IC}{\mathbb{C}}
\newcommand{\IR}{\mathbb{R}}
\newcommand{\IZ}{\mathbb{Z}}

\newcommand{\risomLie}[1]{\vartheta_{#1}}
\newcommand{\risomLieNaive}[1]{\tilde{\vartheta}_{#1}}
\newcommand{\dd}{\mathbf{d}}
\DeclareMathOperator{\ad}{ad}
\newcommand{\adpos}[1]{\operatorname{ad}^+_{#1}}
\newcommand{\adneg}[1]{\operatorname{ad}^-_{#1}}

\newcommand{\rootgr}[1]{U_{#1}}
\newcommand{\rootgrvar}[1]{V_{#1}}
\newcommand{\risomGrp}[1]{\theta_{#1}}

\providecommand\given{}
\newcommand\SetSymbol[1][]{%
	\nonscript\:#1\vert
	\allowbreak
	\nonscript\:
	\mathopen{}}
\DeclarePairedDelimiterX\Set[1]{\{}{\}}{%
	\renewcommand\given{\SetSymbol[\delimsize]}%
	\,#1\, 
}
\DeclarePairedDelimiterX\compactSet[1]{\{}{\}}{%
	\renewcommand\given{\SetSymbol[\delimsize]}%
	#1 
}
\DeclarePairedDelimiterX{\brackets}[1]{(}{)}{#1}
\DeclarePairedDelimiterX{\sqbrackets}[1]{[}{]}{#1}
\DeclarePairedDelimiterX{\abs}[1]{\lvert}{\rvert}{#1}
\DeclarePairedDelimiterX{\gen}[1]{\langle}{\rangle}{\ifblank{#1}{\:\cdot\:}{\renewcommand\given{\SetSymbol[\delimsize]}#1}}

\begin{document}

	\begin{abstract}
		We present \texttt{CubicJordanMatrixAlg}, a GAP package for symbolic computation in cubic Jordan matrix algebras and in related Lie-theoretic structures. As an application, we use it to compute certain (commutator) relations in $F_4$-graded groups that were constructed by De Medts and the author from cubic Jordan matrix algebras.
	\end{abstract}
	
	\maketitle
	
	\section{Introduction}\label{sec:intro}
	
	Probably the most well-known example of a cubic Jordan matrix algebra is
	\begin{equation}\label{eq:herm-def}
		\Her_3(\IC, \IR) \defl \Set*{A \in \brackets*{\begin{smallmatrix}
			\IR & \IC & \IC \\
			\IC & \IR & \IC \\
			\IC & \IC & \IR
		\end{smallmatrix}} \given A_{ij} = \overline{A_{ji}} \text{ for all } i \ne j \in \compactSet{1,2,3}},
	\end{equation}
	the space of Hermitian $(3 \times 3)$-matrices over the complex numbers. On this real vector space, we have a (cubic) norm function $\map{\cubnorm}{\Her_3(\IC,\IR)}{\IR}{}{}$ given by the determinant and a (quadratic) map $ \map{\cubadj}{\Her_3(\IC,\IR)}{\Her_3(\IC,\IR)}{A}{A^\sharp} $ called \emph{adjoint}, which sends a matrix to its adjugate (that is, to the transpose of its cofactor matrix). These two maps equip $\Her_3(\IC, \IR)$ with the structure of a cubic Jordan algebra (which by \cite[34.6]{GPR24} is \enquote{the same} as a cubic norm structure). A similar construction is possible in much higher generality by replacing $\IR$ with a commutative ring $\comring$, $\IC$ with a composition algebra\footnote{There is no universal agreement in the literature on the definition of composition algebras over rings. We follow \cite[19.5]{GPR24}, but the details are not relevant in the context of \cjma.} (such as a quaternion or octonion algebra) $\conic$ over $\comring$ and the complex conjugation with the conjugation of the composition algebra $\conic$. The resulting objects form an important class of examples of so-called \emph{Freudenthal algebras} \cite[39.9~(d)]{GPR24}. Even more generally, we can replace $\IC$ with a multiplicative conic alternative algebra over $\comring$ (a notion we will define below), and the resulting objects $\Her_3(\conic, \comring)$ are called \emph{cubic Jordan matrix algebras}. In this higher generality, the maps $ N $ and $ \sharp $ are defined by explicit formulas which, while straightforward in theory, can be cumbersome to work with in practice. We point out that $ \Her_3(\conic,\comring) $ is not equipped with a multiplication and hence is not an algebra in the usual sense: rather, it is what is called a Jordan algebra \cite[29.1]{GPR24}. In particular, the usual multiplicative structure of matrices will play no role for us.
	
	The purpose of the GAP \cite{GAP4} package \cjma, which we present in this article and which is available at \cite{CJMA}, is the symbolic computation in cubic Jordan matrix algebras. By this we mean that we are not interested in computations in explicit examples such as $\Her_3(\IC, \IR)$, but rather that we want to perform computations in $\Her_3(\conic, \comring)$ for an arbitrary multiplicative conic alternative algebra $\conic$ over an arbitrary commutative ring $\comring$. In other words, we want to perform computations in a \enquote{free} cubic Jordan matrix algebra, but we will not formalise this viewpoint.
	
	The package \cjma{} was born from the need to carry out long computations involving cubic Jordan matrix algebras in \cite{DMW26}, a joint project with Tom De Medts. One of the main goals of \cite{DMW26} is the construction of Lie algebras and groups \enquote{over} arbitrary cubic Jordan matrix algebras and the proof that these objects are graded by the root system $F_4$. We will sketch this application of \cjma{} in \cref{sec:F4}. A brief overview of the ideas underlying \cjma{} has also been given in \cite[Section~9.3]{DMW26}.
	
	Since $\Her_3(\conic, \comring) \cong \comring^3 \oplus \conic^3$ as a $\comring$-module, the problem of (symbolically) computing in $\Her_3(\conic, \comring)$ essentially reduces to the problem of (symbolically) computing in $\conic$ and $\comring$. To better understand the nature of this problem, we briefly consider the related problem of (symbolically) computing in an arbitrary alternative ring (or $ \comring $-algebra) $\alt$ which does not have the additional structure of a multiplicative conic algebra. Recall that a (nonassociative) ring $\alt$ is called \emph{alternative} if it satisfies the alternative laws
	\begin{equation}\label{eq:alternative}
		(aa)b = a(ab) \qquad \text{and} \qquad (ab)b = a(bb)
	\end{equation}
	for all $a,b \in \alt$. Artin's Theorem \cite[Theorem~3.1]{Schafer} says that in any alternative ring, every subring that is generated by at most two elements is associative. This allows us to omit brackets in some expressions, such as $abab$. Further, we have the following Moufang identities \cite[(3.4), (3.5), (3.6)]{Schafer}:
	\begin{equation*}
		(xax)y = x\brackets[\big]{a(xy)}, \qquad y(xax) = \brackets[\big]{(yx)a}x, \qquad (xy)(ax) = x(ya)x
	\end{equation*}
	for all $a,x,y$ in an alternative ring. There are many more non-trivial identities that hold in alternative rings, and their proofs can be rather intricate. Thus it is natural to ask for an algorithm to solve the following problem: Given an element $f$ in the free nonassociative $\comring$-algebra on finitely many generators $x_1, \ldots, x_n$ (so that $f$ is a \enquote{nonassociative noncommutative polynomial in $x_1, \ldots, x_n$} over $\comring$), decide efficiently whether $f$ is an \emph{identity in alternative $ \comring $-algebras}, meaning that $f(a_1, \ldots, a_n) = 0$ for arbitrary $a_1, \ldots, a_n$ in an arbitrary alternative $ \comring $-algebra. This problem is unsolved in this generality, but we report on some prior work in \cref{sec:literature}.
	
	Returning to \texttt{CubicJordanMatrixAlg}, we first have to clarify the pivotal notion of multiplicative conic alternative algebras. We follow the notation and exposition in \cite[Sections~16,~17]{GPR24}.
	
	\begin{convention}
		Throughout, rings and algebras are assumed to be unital but not necessarily associative or commutative. However, by \enquote{commutative ring} we always mean \enquote{associative commutative ring}.
	\end{convention}
	
	Recall that a \emph{quadratic form} on a module $ M $ over a commutative ring $ \comring $ is a map $ \map{q}{M}{\comring}{}{} $ such that $ q(tv) = t^2 q(v) $ for all $ t \in \comring $, $ v \in M $ and the map $ \map{Dq}{M \times M}{\comring}{(v,w)}{q(v+w)-q(v)-q(w)} $ is $ \comring $-bilinear. We call $ Dq $ the \emph{linearisation of $ q $}, and will always write $ q(v,w) $ for $ Dq(v,w) $.
	
	\begin{definition}[{\cite[16.1]{GPR24}}]
		A \emph{conic $\comring$-algebra} is a pair $(\conic, \connorm)$ consisting of a $\comring$-algebra $\conic$ and a quadratic form $\map{\connorm}{\conic}{\comring}{}{}$ called the \emph{norm} such that $\connorm(1_\conic) = 1_\comring$ and the \emph{Cayley-Hamilton equation}
		\[ a^2 - \connorm(1_\conic,a)a + \connorm(a)1_\conic = 0 \]
		holds for all $a \in \conic$. The map $\map{\contr}{\conic}{\comring}{a}{\connorm(1_\conic,a)}$ is called the \emph{trace of $(\conic,\connorm)$} and $\map{\conconj{}}{\conic}{\conic}{a}{\contr(a)1_\conic -a}$ is called the \emph{conjugation of $(\conic,\connorm)$}. The conic algebra $(\conic,\connorm)$ is called \emph{multiplicative} if $\connorm(ab) = \connorm(a) \connorm(b)$ for all $a,b \in \conic$ and it is called \emph{alternative} if it satisfies~\eqref{eq:alternative} for all $a,b \in \conic$.
	\end{definition}
	
	\begin{convention}
		For the rest of this article, we will simply write \enquote{conic algebra} for \enquote{multiplicative conic alternative algebra}. We will also refer to $\conic$ as a conic algebra, leaving the norm $\connorm$ implicit.
	\end{convention}
	
	\begin{example}
		Composition algebras over commutative rings (such as the classical examples of complex numbers, quaternions and octonions over the real numbers) are conic algebras \cite[19.13]{GPR24}. For any commutative ring $ \comring $, any $ \comring $-algebra that is projective of rank~$ 2 $ is a conic algebra \cite[16.4]{GPR24}.
	\end{example}
	
	In addition to the intricate identities that are satisfied in alternative rings, any conic algebra $\conic$ satisfies identities involving its trace, norm and conjugation, such as the following.
	\begin{miscthm}[Identities]\label{conic-identities}
		For all $a,b,c$ in a conic algebra $ \conic $ over a commutative ring $ \comring $, the following hold:
		\begin{enumerate}[(i)]
			\item $\conconj{a+b}=\conconj{a}+\conconj{b}$, $\conconj{ab} = \conconj{b} \conconj{a}$, $\conconj{\conconj{a}}=a$ and $\conconj{ta} = t \conconj{a}$ for $t \in \comring$. \cite[17.2~(a), (16.5.4)]{GPR24}
			
			\item\label{id:tr} $\contr(ab) = \contr(ba)$, $\contr\brackets[\big]{(ab)c} = \contr\brackets[\big]{a(bc)}$, $\contr(\conconj{a}) = \contr(a)$. \cite[(16.5.7), 16.13, 17.12~(a)]{GPR24}
			
			\item\label{id:n} $\connorm(\conconj{a}) = \connorm(a)$, $\connorm(a,b)=\contr(a\conconj{b})$, $\connorm(1,a) = \contr(a)$. \cite[(16.5.2), (16.5.7), (16.12.3), 17.12~(a)]{GPR24}
			
			\item\label{id:conj} $\contr(a)1_\conic = a + \conconj{a}$ and $\connorm(a)1_\conic = a \conconj{a} = \conconj{a}a$. \cite[(16.5.6)]{GPR24}
			
			
			\item\label{id:kirmse} Kirmse's identity: $a \connorm(b) = (a \conconj{b})b$. \cite[(17.4.1)]{GPR24}
			
			\item\label{id:tr-square} $\contr(a^2) = \contr(a)^2 - 2\connorm(a)$. \cite[(16.5.8)]{GPR24}
		\end{enumerate}
		Note that some of these are equations in $\conic$ while others are equations in the base ring $\comring$, but we refer to both kinds of relations as \emph{identities in conic algebras}.
	\end{miscthm}

	\begin{miscthm}[Goal]
		As for alternative rings, we would like to decide whether a given expression is an identity in conic algebras---that is, if it is valid in every conic algebra (or in the ground ring of any conic algebra). The purpose of \texttt{CubicJordanMatrixAlg} is more humble: Given an expression $f$ for which we would like to decide whether it is an identity in conic algebras, \texttt{CubicJordanMatrixAlg} can compute a \enquote{simplified} expression $\tilde{f}$ such that $f$ is an identity in conic algebras if and only if $\tilde{f}$ is. In many examples that arise in applications, we have $\tilde{f}=0$, in which case \texttt{CubicJordanMatrixAlg} has proven that $f$ is an identity. In other cases, it remains to prove by hand that $\tilde{f}$ is an identity, but this is usually much easier than for $f$. However, \texttt{CubicJordanMatrixAlg} cannot prove that an expression is \emph{not} an identity.
	\end{miscthm}

	\section{A first usage example}\label{sec:ex}
	
	Let us illustrate how to use \texttt{CubicJordanMatrixAlg} for computations in conic algebras. After following the installation instructions in \cite{CJMA}, the package should be loaded in GAP and initialised with \texttt{InitCJMA}.
\begin{verbatim}
gap> LoadPackage("CubicJordanMatrixAlg");
gap> InitCJMA(1, 2, 3, true);
\end{verbatim}
	We will explain the meaning of the arguments $ m_1,m_2,m_3 $ (positive integers), \texttt{user\_vars} (boolean) of \texttt{InitCJMA} throughout this section. The first two arguments specify the maximal number of pairwise independent elements \texttt{t1}, \dots, \texttt{t}$ m_1 $ in $ \comring $ and \texttt{a1}, \dots, \texttt{a}$ m_2 $ in $ \conic $ that we want to work with. We think of these elements as \enquote{indeterminates in $ \comring $ and $ \conic $}, and as we will see in \cref{sec:implement}, this is indeed how they are implemented internally. In our example, indeterminates \texttt{t1} in $ \comring $ and \texttt{a1}, \texttt{a2} in $ \conic $ are defined. We may access the $ i $-th indeterminates using \texttt{ComRingIndet(i)} and \texttt{ConicAlgIndet(i)}. If the fourth argument \texttt{user\_vars} passed to \texttt{InitCJMA} is \texttt{true}, then GAP variables with the same names as the indeterminates are defined for the user's convenience. In our example, this is equivalent to the following definitions:
\begin{verbatim}
gap> t1 := ComRingIndet(1);; a1 := ConicAlgIndet(1);; a2 := ConicAlgIndet(2);;
\end{verbatim}

	GAP functions to (symbolically) compute the conjugation, norm and trace, which make use of some of the relations in~\ref{conic-identities}, are available:
\begin{verbatim}
gap> ConicConj(a1+a1*a2); ConicTr(a1+a1*a2); ConicNorm(a1*a2); ConicNorm(a1+a2);
(1)*a1'+(1)*(a2'*a1')
tr(a1)+tr(a1a2)
n(a1)*n(a2)
n(a1)+n(a2)+tr(a1a2')
\end{verbatim}
	Here \texttt{a1'}, \texttt{a2'} are new symbols representing the conjugates of \texttt{a1}, \texttt{a2} while \texttt{tr(...)} and \texttt{n(...)} are new symbols representing the trace or norm of their argument, respectively. Note that the multiplicativity of $ \connorm $ allows us to work with only $ m_2 $ many symbols \texttt{n(a1)}, \dots, \texttt{n(a}$ m_2 $\texttt{)} to represent the norm. The trace $ \contr $, however, is not multiplicative, so we also need symbols such as \texttt{tr(a1a2)}. In theory we would like to have a symbol \texttt{tr(b)} for each of the infinitely many \enquote{monomials} \texttt{b} in $ \conic $, but for implementation reasons, the number of symbols must be bounded and known during initialisation. Because of this, the third argument $ m_3 $ of \texttt{InitCJMA} specifies the maximal length of monomials \texttt{b} that are allowed in a symbol \texttt{tr(b)}. In our example with $ m_3=3 $, we can compute \texttt{ConicTr(a1*a1*a2)}, but \texttt{ConicTr(a1*a1*a2*a2)} throws an error.
	
	The values $ m_1 $, $ m_2 $, $ m_3 $ significantly affect the time needed for initialisation and, to a small extent, also the runtime of \cjma. The effect of $ m_2 $ and $ m_3 $ is particularly strong because the number of required indeterminates \texttt{tr(b)} grows quickly. 
	
	The expressions computed by \cjma{} are not simplified unless the user asks for it. This is done with \texttt{Simplify}, which is the main function of \texttt{CubicJordanMatrixAlg}.
\begin{verbatim}
gap> x := a1*a2+a1*ConicConj(a2)-ConicTr(a2)*a1;
(-tr(a2))*a1+(1)*(a1*a2)+(1)*(a1*a2')
gap> Simplify(x);
0_C
\end{verbatim}
	Indeed, the term $x \defl a_1 a_2 + a_1 \conconj{a_2} - \contr(a_2)a_1$ is zero for arbitrary $a_1$, $a_2$ in any conic algebra by~\ref{conic-identities}~\ref{id:conj}. \texttt{Simplify} works by analysing the input \texttt{x} for ways in which (a fixed list of) identities such as those in~\ref{conic-identities} can be used to \enquote{reduce} \texttt{x}. In fact, the identity $\contr(a_1)1_\conic = a_1 + \conconj{a_1}$ is one of the identities used by \texttt{Simplify}, so the example above is a somewhat trivial application. The true power of \texttt{Simplify} lies in the fact that it can be used for much longer expressions \texttt{x}, and that we can use it within the programming environment of GAP to analyse a large number of terms \texttt{x}. For a non-trivial application, see \cref{sec:F4} and \cite{DMW26}.

%
%
%
%
%
%
%
%
%
	
	\section{Cubic Jordan matrix algebras}\label{sec:cjma}
	
	\begin{notation}
		For the rest of this article, we denote by $\comring$ a commutative ring and by $(\conic, \connorm)$ a conic algebra with trace $\contr$ and conjugation $\conconj{}$. We also fix a triple $\Gamma=(\gamma_1, \gamma_2, \gamma_3) \in \comring^3$ of structure constants and we denote by $\perm \defl \Set{(1,2,3),(2,3,1), (3,1,2)}$ the set of cyclic permutations of $(1,2,3)$. For simplicity, we assume that $ \gamma_1,\gamma_2,\gamma_3 $ are invertible, though some of our statements remain valid without this assumption.
	\end{notation}
	
	In this section, we briefly define the notion of cubic Jordan matrix algebra as in \cite[36.4]{GPR24}. They provide examples of cubic Jordan algebras. The precise axioms of the latter are irrelevant for our purposes, but can be found in \cite[33.1, 33.4]{GPR24} and \cite[2.4]{DMW26} (for cubic norm structures, which are \enquote{equivalent to} cubic Jordan algebras by \cite[34.6]{GPR24}) and \cite[29.1]{GPR24} (for general Jordan algebras). We will also show how to perform symbolic computations in cubic Jordan matrix algebras with \cjma.
	
	\begin{remark}\label{matrix-shape}
		Denote by $ e_{ij} $ the $ (3 \times 3) $-matrix with $ 1 $ at position $ (i,j) $ and $ 0 $ everywhere else. Put $ e_i \defl e_{ii} $ and $ \cubel{u}{jl} \defl \gamma_l u e_{jl} + \gamma_j \conconj{u} e_{lj} $ for all $ u \in \conic $, $ (i,j,l) \in \perm $. If the diagonal matrix $ \Gamma= \operatorname{diag}(\gamma_1, \gamma_2, \gamma_3) $ is invertible and $ \conic $ is a faithful $ \comring $-algebra (so that we may identify $ \comring $ with $ \comring 1_\conic \subseteq \conic $), then we may define $ \Her_3(\conic, \comring, \Gamma) $ as the $ \comring $-module of $ (3 \times 3) $-matrices over $ \conic $ that are fixed under (the involution) $ \map{}{}{}{A}{\Gamma^{-1} \conconj{A}^T \Gamma} $ and whose diagonal elements lie in $\comring 1_\conic$ \cite[36.2]{GPR24}. Any such matrix may be written in a unique way as $ \xi_1 e_1 + \xi_2 e_2 + \xi_3 e_3 + \cubel{u_1}{23} + \cubel{u_2}{31} + \cubel{u_3}{12} $ for $ \xi_1, \xi_2, \xi_3 \in \comring $, $ u_1, u_2, u_3 \in \conic $. This observation motivates the following more general definition of $ \Her_3(\conic, \comring, \Gamma) $, which is equivalent to the previous one if $ \Gamma $ is invertible and $ \conic $ is faithful.
%
	\end{remark}
	
	\begin{definition}
		As a $\comring$-module, the \emph{cubic Jordan matrix algebra over $\conic$ with respect to $\Gamma$} is $\Her_3(\conic, \comring, \Gamma) \defl \comring^3 \oplus \conic^3$. We will write an element $ (\xi_1, \xi_2, \xi_3, u_1, u_2, u_3) $ of $ \Her_3(\conic, \Gamma) $ (with $ \xi_1,\xi_2, \xi_3 \in \comring $ and $ u_1,u_2,u_3 \in \conic $) as a a formal expression
		\[ \xi_1 e_1 + \xi_2 e_2 + \xi_3 e_3 + \cubel{u_1}{23} + \cubel{u_2}{31} + \cubel{u_3}{12} = \sum_{(i,j,l) \in \perm} \brackets[\big]{\xi_i e_i + \cubel{u_i}{jl}}. \]
		We also put $ \cubel{a}{ji} \defl \cubel{\conconj{a}}{ij} $ for $ a \in \conic $ and $ (i,j,l) \in \perm $ as well as $\cubel{t}{ii} \defl te_i$ for $t \in \comring$ and $i \in \Set{1,2,3}$. If $ \Gamma = (1,1,1) $, we simply write $ \Her_3(\conic, \comring) $ for $ \Her_3(\conic, \comring, \Gamma) $, as in \cref{sec:intro}. Further, we will usually write $ \Her_3(\conic, \Gamma) $ for $ \Her_3(\conic, \comring, \Gamma) $, leaving the base ring $ \comring $ implicit. 
	\end{definition}
	
	
	In the following, we put $ \jor \defl \Her_3(\conic, \Gamma) $. As a $\comring$-module, $J$ does not depend on $\Gamma$, but we define maps
	\begin{align*}
		\map{\cubnorm}{\jor}{\comring}{}{}, \quad \map{\cubadj}{\jor}{\jor}{}{}, \quad \map{\cubtr}{\jor \times \jor}{\comring}{}{}, \quad \map{\cubtimes}{\jor \times \jor}{\jor}{}{}
	\end{align*}
	that do depend on $ \Gamma $ and which
	turn $ \jor $ into a cubic norm structure. They are given by the following explicit formulas in \cite[36.4]{GPR24}:
	\begin{align*}
		\cubnorm(x) &\defl \xi_1 \xi_2 \xi_3 + \gamma_1 \gamma_2 \gamma_3 \contr(u_1 u_2 u_3) - \sum_{(i,j,l) \in \perm} \gamma_j \gamma_l t_i \connorm(u_i), \\ 
		x^\sharp &\coloneq \sum_{(i,j,l) \in \perm} \brackets[\Big]{\brackets[\big]{\xi_j \xi_l - \gamma_j \gamma_l \connorm(u_i)} e_i + \cubel{(-\xi_i u_i + \gamma_i \conconj{(u_j u_l)})}{jl}}, \\
		x \times y &\coloneq \sum_{(i,j,l) \in \perm} \Bigl(\brackets[\big]{\xi_j \eta_l + \eta_j \xi_l - \gamma_j \gamma_l \connorm(u_i, v_i)} e_i + \cubel{\brackets[\big]{-\xi_i v_i - \eta_i u_i + \gamma_i \conconj{(u_j v_l + v_j u_l)}}}{jl}\Bigr), \\
		T(x,y) &\coloneq \sum_{(i,j,l) \in \perm} \brackets[\big]{\xi_i \eta_i + \gamma_j \gamma_l \connorm(u_i, v_i)}
	\end{align*}
	for all $ x = \sum_{(i,j,l) \in \perm} (\xi_i e_i + \cubel{u_i}{jl}) \in \jor, y = \sum_{(i,j,l) \in \perm} (\eta_i e_i + \cubel{v_i}{jl}) \in \jor $. The maps $ N $, $ \sharp $, $ T $ are called the norm, adjoint and (bilinear) trace of $ J $, respectively. In the case $ \comring = \IR $, $ \conic=\IC $, and $ \Gamma = (1,1,1) $, they correspond precisely to the determinant, the adjugate map and the Killing form $ (x,y) \mapsto \operatorname{Tr}(xy) $ where $ \operatorname{Tr} $ denotes the usual matrix trace.
	
	By \cite[33.9]{GPR24}, the map
	\begin{align*}
		\map{U}{\jor \times \jor}{\jor&}{(x,y)}{U_x y \defl \cubtr(x,y)x - x^\cubadj \cubtimes y}
	\end{align*}
	gives $ \jor $ the structure of a Jordan algebra 
	(in the sense of \cite[29.1]{GPR24}). Its linearisation is
	\begin{align*}
		\map{\tribrace{}{}{}}{\jor \times \jor \times \jor}{\jor&}{(x,y,z)}{U_{x+z}y - U_x y - U_z y = \cubtr(x,y)z + \cubtr(y,z)x - (z \cubtimes x) \cubtimes y}.
	\end{align*}
	We will also write $ D_{x,y}(z) \defl \tribrace{x}{y}{z} $ for $ x,y,z \in J $.
	
	\begin{example}[Peirce decomposition]\label{ex:peirce}
		For all $ (i,j,l) \in \perm $, we denote by $ J_{ii} \defl \comring e_i $ and $ J_{ij} \defl J_{ji} \defl \cubel{\conic}{ij} $ the components of $ \jor $. These are the Peirce components of $ \jor $ by \cite[37.8]{GPR24}, so it follows from the general theory of Peirce decompositions in Jordan algebras \cite[32.15]{GPR24} that
		\begin{equation}\label{eq:peirce-decomp}
			\tribrace{\jor_{ij}}{\jor_{jl}}{\jor_{lm}} \subseteq J_{im} \qquad \text{for all }  i,j,l,m \in \Set{1,2,3}.
		\end{equation}
		Using that $\jor_{ij} = \jor_{ji}$ for all $ i,j \in \Set{1,2,3} $, we can also apply~\eqref{eq:peirce-decomp} in situations where it is not immediately apparent. For example, $\tribrace{\jor_{12}}{\jor_{12}}{\jor_{32}} \subseteq \jor_{23}$.
		Further, one can show that $\tribrace{\jor_{ij}}{\jor_{lm}}{\jor_{pq}} = \compactSet{0}$ if it is not possible to rearrange the indices $i,j,l,m,p,q \in \Set{1,2,3}$ in a way which allows us to apply~\eqref{eq:peirce-decomp}.
		
		Alternatively, we can of course use the explicit nature of the formulas defining $ \cubtr $ and $ \cubtimes $ (and hence $ \tribrace{}{}{} $) to prove the above statements by means of a direct computation. Let us illustrate how to do this with \cjma. For this we need a higher number of indeterminates than in \cref{sec:ex}, so let us start a new GAP session and initialise the package accordingly.
\begin{verbatim}
gap> LoadPackage("CubicJordanMatrixAlg");;
gap> InitCJMA(6, 3, 4, true);
\end{verbatim}
		We can create elements of $ \jor $ using the function \texttt{CubicEl}. Note that \cjma{} displays $ t_1 e_2 $ as $ \cubel{t_1}{22} $ and that $ \cubel{a_2}{21} $ is converted to the equivalent expression $ \cubel{\conconj{a_2}}{12} $.
\begin{verbatim}
gap> CubicEl(t1, 2, 2); CubicEl(a2, 2, 1);
(t1)[22]
((1)*a2')[12]
\end{verbatim}
		The $ D $-operation $ D_{x,y}(z) = \tribrace{x}{y}{z} $ can be computed with \texttt{JordanD}. Here \texttt{g1}, \texttt{g2}, \texttt{g3} stand for the fixed (but arbitrary) structure constants $ \gamma_1, \gamma_2, \gamma_3 \in \comring $.
\begin{verbatim}
gap> JordanD(CubicEl(t1, 2, 2), CubicEl(a2, 2, 1), CubicEl(a3, 1, 2));
(g1*g2*t1*tr(a2a3))[22]
\end{verbatim}
		This computation shows not only that $ \tribrace{\jor_{22}}{\jor_{21}}{\jor_{12}} \subseteq J_{22} $, but in fact proves the stronger formula $ \tribrace{\cubel{t_1}{22}}{\cubel{a_2}{21}}{\cubel{a_3}{12}} = \cubel{\gamma_1 \gamma_2 t_1 \contr(a_2 a_3)}{22} $ for $ a_1, a_2, a_3 \in \conic $. Using \cjma, we have in \cite[9.15]{DMW26} computed similar explicit formulas for all cases that may arise in \eqref{eq:peirce-decomp}. We point out that, while these computations could still be done by hand with reasonable effort and patience, we are not aware of a reference which proves these formulas for cubic Jordan matrix algebras (with the exception of the single formula in \cite[(37.7.8)]{GPR24}). A noteworthy reference in this context is \cite[(24)]{McCrimmon66} which on first glance contains the same formulas as in \cite[9.15]{DMW26} for the special case $ \gamma_1=\gamma_2=\gamma_3=1 $. However, the setup is actually different: The formulas in \cite[(24)]{McCrimmon66} are used to \emph{define} a Jordan algebra over an alternative ring with involution, but this does not imply that cubic Jordan matrix algebras in our sense satisfy the same formulas.
		
		We now prove the formula $ \tribrace{\cubel{a_1}{23}}{\cubel{a_2}{31}}{\cubel{a_3}{13}} = \cubel{\gamma_1 \gamma_3 (a_1 a_2) a_3}{23} $ from \cite[9.15~(v)]{DMW26} as an example of a computation in which some final steps by hand are necessary. Indeed, \cjma{} yields only the following:
\begin{verbatim}
gap> JordanD(CubicEl(a1, 2, 3), CubicEl(a2, 3, 1), CubicEl(a3, 1, 3));
((g1*g3*tr(a2a3))*a1+(-g1*g3)*((a1*a3')*a2'))[23]
\end{verbatim}
		Kirmse's identity (\ref{conic-identities}~\ref{id:kirmse}) says that $ a_1 \connorm(b) = (a_1 \conconj{b}) b $ for all $ a_1, b \in \conic $. We can linearise this identity: We put $ b \defl a_2+a_3 $ to obtain
		\begin{align*}
			a_1 \connorm(a_2) + a_1 \connorm(a_3) + a_1 \connorm(a_2, a_3) = (a_1 \conconj{a_2}) a_2 + (a_1 \conconj{a_3}) a_3 + (a_1 \conconj{a_2}) a_3 + (a_1 \conconj{a_3}) a_2
		\end{align*}
		and apply Kirmse's identity two more times to find that $ a_1\connorm(a_2,a_3) = (a_1 \conconj{a_2}) a_3 + (a_1 \conconj{a_3}) a_2 $. Here $ \connorm(a_2,a_3) = \contr(a_2 \conconj{a_3}) = \contr(\conconj{a_2} a_3) $ by~\ref{conic-identities}~\ref{id:n} and~\ref{conic-identities}~\ref{id:tr}. Replacing $ a_2 $ by $ \conconj{a_2} $, we thus have $ a_1 \contr(a_2 a_3) = (a_1 a_2) a_3 + (a_1 \conconj{a_3}) \conconj{a_2} $. This implies that the expression for $ \tribrace{\cubel{a_1}{23}}{\cubel{a_2}{31}}{\cubel{a_3}{13}} $ which we computed with \cjma{} is indeed equal to $ \cubel{\gamma_1 \gamma_3 (a_1 a_2) a_3}{23} $.
	\end{example}

	\section{\texorpdfstring{$F_4$}{F4}-graded Lie algebras and groups}\label{sec:F4}
	
	One of the main results of \cite{DMW26} is the construction of a root graded Lie algebra $ L $ of type $ F_4 $ from an arbitrary conic algebra $ \conic $ and the construction of a root graded group $ G $ of type $F_4$ inside $ \Aut(L) $ that is \emph{coordinatised by $ \conic $}. (We will define the notion of root gradings in \ref{def:rg-lie},~\ref{rgg} and the notion of coordinatisations in~\ref{rgg},~\ref{ex-prob}.)
	\cjma{} supports computations in both $ L $ and $ G $ and is used for this purpose in \cite{DMW26}. In this section, we give a rough sketch of this application of \cjma{}.
	
	We first have to recall the notion of root systems from Lie theory.
	
	\begin{definition}
		A \emph{root system} is a finite subset $\roots$ of a Euclidean space $V$ such that $0 \notin \roots$ and for all $\alpha \in \roots$, we have $\roots^{\reflbr{\alpha}} = \roots$ where $\map{\reflbr{\alpha}}{V}{V}{\beta}{\beta - 2 \frac{\alpha \cdot \beta}{\alpha \cdot \alpha} \alpha}$ is the reflection along the hyperplane $\alpha^\perp$. For any root system $\roots$, we put $\roots^0 \defl \roots \cup \compactSet{0}$.
	\end{definition}
	
\begin{figure}[htb]
		\newcommand{\smallboxlength}{9mm}%
		\newcommand{\bigboxlength}{22mm}%
		\centering$ \scalebox{.7}{%
			\begin{tikzpicture}[x=36mm, y=25mm, 	label distance=-3pt]
				\node[mysquare=\smallboxlength] at (0,3) (N03) {$\bar{2}000$};
				\node[mysquare=\smallboxlength] at (1,1) (N11) {$\bar{1}111$};
				\node[mysquare=22mm] at (1,2) (N12) {$\begin{matrix} \bar{1}\bar{1}11 & \bar{1}001 & \bar{1}010 \\[.5ex] & \bar{1}1\bar{1}1 & \bar{1}100 \\[.5ex] && \bar{1}11\bar{1} \end{matrix}$};
				\node[mysquare=\bigboxlength] at (1,3) (N13) {$\begin{matrix} \bar{1}1\bar{1}\bar{1} & \bar{1}00\bar{1} & \bar{1}0\bar{1}0 \\[.5ex] & \bar{1}\bar{1}1\bar{1} & \bar{1}\bar{1}00 \\[.5ex] && \bar{1}\bar{1}\bar{1}1 \end{matrix}$};
				\node[mysquare=\smallboxlength] at (1,4) (N14) {$\bar{1}\bar{1}\bar{1}\bar{1}$};
				\node[mysquare=\bigboxlength] at (2,1) (N21) {$\begin{matrix} 0200 & 0110 & 0101 \\[.5ex] & 0020 & 0011 \\[.5ex] && 0002 \end{matrix}$};
				\node[mysquare=\bigboxlength] at (2,2) (N22) {$\begin{matrix} & 01\bar{1}0 & 010\bar{1} \\[.5ex] 0\bar{1}10 & 0000 & 001\bar{1} \\[.5ex] 0\bar{1}01 & 00\bar{1}1 & \end{matrix}$};
				\node[mysquare=\bigboxlength] at (2,3) (N23) {$\begin{matrix} 0\bar{2}00 & 0\bar{1}\bar{1}0 & 0\bar{1}0\bar{1} \\[.5ex] & 00\bar{2}0 & 00\bar{1}\bar{1} \\[.5ex] && 000\bar{2} \end{matrix}$};
				\node[mysquare=\smallboxlength] at (3,0) (N30) {$1111$};
				\node[mysquare=\bigboxlength] at (3,1) (N31) {$\begin{matrix} 1\bar{1}11 & 1001 & 1010 \\[.5ex] & 11\bar{1}1 & 1100 \\[.5ex] && 111\bar{1} \end{matrix}$};
				\node[mysquare=\bigboxlength] at (3,2) (N32) {$\begin{matrix} 11\bar{1}\bar{1} & 100\bar{1} & 10\bar{1}0 \\[.5ex] & 1\bar{1}1\bar{1} & 1\bar{1}00 \\[.5ex] && 1\bar{1}\bar{1}1 \end{matrix}$};
				\node[mysquare=\smallboxlength] at (3,3) (N33) {$1\bar{1}\bar{1}\bar{1}$};
				\node[mysquare=\smallboxlength] at (4,1) (N41) {$2000$};
				\path[ugentred]
					(0,4.5) node (C0) {}
					(1,4.5) node (C1) {}
					(2,4.5) node (C2) {}
					(3,4.5) node (C3) {}
					(4,4.5) node (C4) {};
				\path[ugentblue]
					(4.5,0) node (R0) {}
					(4.5,1) node (R1) {}
					(4.5,2) node (R2) {}
					(4.5,3) node (R3) {}
					(4.5,4) node (R4) {};
				\draw[myedge,ugentblue]
					(N11) -- (N21) -- (N31) -- (N41)
					(N12) -- (N22) -- (N32)
					(N03) -- (N13) -- (N23) -- (N33);
				\draw[myedge,ugentblue,opacity=.15]
					(-.25,0) -- (N30) -- (R0)
					(-.25,1) -- (N11) (N41) -- (R1)
					(-.25,2) -- (N12) (N32) -- (R2)
					(-.25,3) -- (N03) (N33) -- (R3)
					(-.25,4) -- (N14) -- (R4);
				\draw[myedge,ugentred]
					(N11) -- (N12) -- (N13) -- (N14)
					(N21) -- (N22) -- (N23)
					(N30) -- (N31) -- (N32) -- (N33);
				\draw[myedge,ugentred,opacity=.15]
					(0,-.25) -- (N03) -- (C0)
					(1,-.25) -- (N11) (N14) -- (C1)
					(2,-.25) -- (N21) (N23) -- (C2)
					(3,-.25) -- (N30) (N33) -- (C3)
					(4,-.25) -- (N41) -- (C4);
				\foreach \x/\val in {0/-2, 1/-1, 2/0, 3/1, 4/2}{
					\node[ugentred, opacity=1] at (\x, -0.4){$\val$};
				}
				\foreach \y/\val in {0/2, 1/1, 2/0, 3/-1, 4/-2}{
					\node[ugentblue, opacity=1] at (-0.4, \y){$\val$};
				}
			\end{tikzpicture}
		} $
		\caption{Each rectangle represents one of the elements of $G_2^0 \subseteq \IZ^2$. Big rectangles correspond to short roots and small rectangles to long roots. Each four-digit number $\epsilon_1 \epsilon_2 \epsilon_3 \epsilon_4$ represents the element $\sum_{i=1}^4 \epsilon_i e_i$ of $F_4^0$ where $\bar{2}$ and $\bar{1}$ stand for $-2$ and $-1$, respectively. For any $\alpha \in F_4^0$, $\pi(\alpha)$ corresponds to the rectangle in which $\alpha$ lies.}
		\label{fig:F4}
	\end{figure}
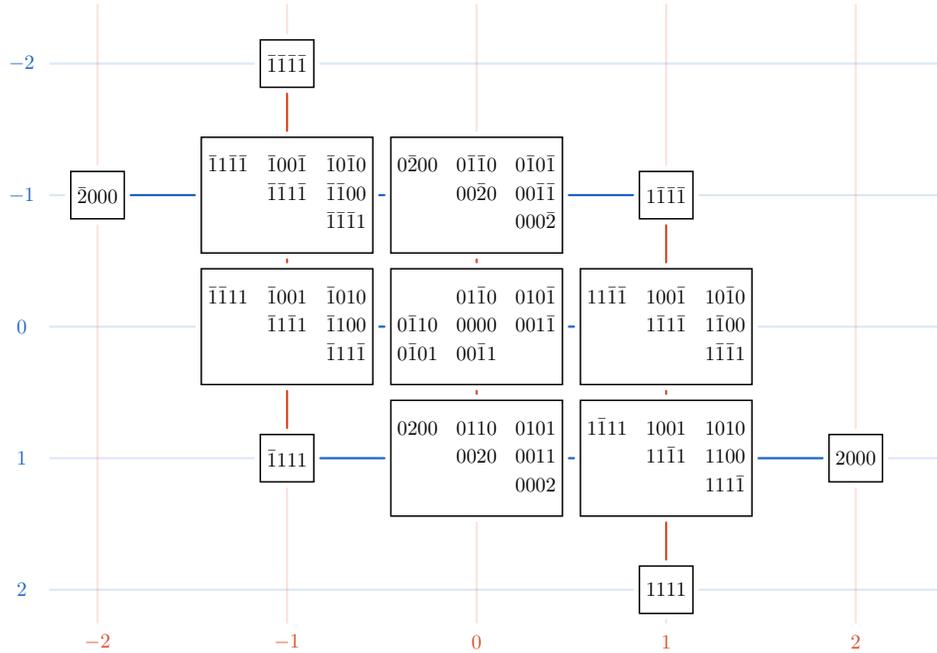

	\begin{miscthm}[The root systems $ G_2 $ and $F_4$]\label{G2F4}
		We consider the root system $G_2$ in the following representation as a subset of the Euclidean space $ \IR^2 $ with the inner product given by the Gram matrix $ \brackets*{\begin{smallmatrix}
			2 & -1 \\
			-1 & 2
		\end{smallmatrix}} $:
		\begin{align*}
			G_2 = \Set[\big]{\pm (a,b) \given (a,b) \in \compactSet{(-2,-1), (-1,-1), (-1,-2), (0,-1), (1,-1), (1,0)}}.
		\end{align*}
		The roots $ \pm (-1,-1) $, $ \pm (0,-1) $, $ \pm (1,0) $ are short (of length $ \sqrt{2} $) while the others are long (of length $ \sqrt{6} $).
		
		Similarly, we use the following representation of the root system $F_4$ as a subset of $\IR^4$ with the standard inner product, whose standard basis we denote by $(e_1, \ldots, e_4)$:
		\begin{align*}
			F_4 &\defl \Set{\pm 2e_i \given i \in \compactSet{1,\ldots, 4}} \cup \Set*{\sum_{i=1}^4 \epsilon_i e_i \given \epsilon_1, \ldots, \epsilon_4 \in \compactSet{\pm 1}} \\
			&\qquad \mathord{}\cup \Set{a e_i + b e_j \given i \ne j \in \compactSet{1, \ldots, 4}, a,b \in \compactSet{\pm 1}}
		\end{align*}
		Here the roots in the first row are the long ones while those in the second row are the short ones.
		
		Both root systems have \emph{5-gradings}: Decompositions $G_2 = \bigsqcup_{i=-2}^2 \Ggrad{i}$ and $F_4 = \bigsqcup_{i=-2}^2 \Fgrad{i}$ defined by
		\begin{align*}
			\Ggrad{i} \defl \Set{(a,b) \in G_2 \given a=i}, \qquad \Fgrad{i} \defl \Set{(a_1,\ldots, a_4) \in F_4 \given a_1=i}.
		\end{align*}
		
		We can define a linear surjection $ \map{\pi}{F_4^0}{G_2^0}{(p,i,j,l)}{\brackets[\big]{p, (p+i+j+l)/2}} $ which relates the two root systems. This map is connected to a certain Tits index, which is studied in more detail in \cite[Section~8]{DMW26}. The values of $\pi$ are displayed graphically in \cref{fig:F4}.
	\end{miscthm}

	\premidfigure
	
	\postmidfigure
	
	We now turn to (computations in) the Lie algebra $L$ from \cite{DMW26} that can be constructed from $J$.
	
	\begin{definition}\label{def:rg-lie}
		Let $ L $ be a Lie algebra over $ \comring $ and let $ \roots $ be a finite root system. A \emph{$ \roots $-grading of $ L $} is a family $ (L_\alpha)_{\alpha \in \roots^0} $ of submodules such that $ L = \bigoplus_{\alpha \in \roots^0} L_\alpha $ as $ \comring $-modules and such that for all $ \alpha,\beta \in \roots^0 $, we have $ [L_\alpha, L_\beta] \subseteq L_{\alpha+\beta} $ if $ \alpha+\beta \in \roots^0 $ and $ [L_\alpha,L_\beta] = \compactSet{0} $ otherwise.
	\end{definition}
	
	\begin{miscthm}[An $ F_4 $-graded Lie algebra]\label{lie-alg}
		We denote by $L$ the Lie algebra that is constructed from $J$ in \cite[Section~4]{DMW26}. This Lie algebra has a $ G_2 $-grading $ (L_\alpha)_{\alpha \in G_2^0} $ and an $ F_4 $-grading $ (L_\alpha)_{\alpha \in F_4^0} $. For $ (i,j) \in G_2^0 $, we will write $ L_{ij} $ or $L_{i,j} $ for $ L_{(i,j)} $. We have isomorphisms of $ \comring $-modules
		\begin{align*}
			\map{\risomLieNaive{\alpha}&}{\comring}{L_\alpha}{}{} \quad \text{for } \alpha \in G_2 \text{ long}, & \map{\risomLieNaive{\beta}&}{\jor}{L_\alpha}{}{} \quad \text{for } \alpha \in G_2 \text{ short,} \\
			\map{\risomLieNaive{\alpha}&}{\comring}{L_\beta}{}{} \quad \text{for } \beta \in F_4 \text{ long}, & \map{\risomLieNaive{\beta}&}{\conic}{L_\beta}{}{} \quad \text{for } \beta \in F_4 \text{ short.}
		\end{align*}
		The $ 0 $-spaces $ L_{00} $ and $ L_{0000} $ are more difficult to describe. As $ \comring $-modules, they are generated by certain elements $ \dd_{c,c'} = \dd(c,c') $ that are defined for all $ c,c' \in \jor $ and by additional elements $ \xi $, $ \zeta $.
		The map $ (c,c') \mapsto \dd_{c,c'}$ is $\comring$-bilinear.
		The two gradings of $ L $ are related by the fact that $ L_\beta = \bigoplus_{\alpha \in \pi^{-1}(\beta)} L_\alpha $ for all $ \beta \in G_2^0 $.
		For any short $ \beta \in G_2 $, we have $ \abs{\pi^{-1}(\beta)} = 6 $ and the decomposition of $ L_\beta $ into $ F_4 $-root spaces corresponds precisely to the natural decomposition $ \jor = \comring^3 \oplus \conic^3 $. This means that for $\alpha \in F_4$, the $\comring$-module $L_\alpha$ is one of the six canonical direct summands of $L_{\pi(\alpha)} \cong J$. In \cref{fig:F4}, the placement of the $ F_4 $-roots in $\pi^{-1}(\beta)$ for short $\beta \in G_2$ mimics the (informal) matrix shape of $J$ (see \cref{matrix-shape}):
		For example, $\risomLieNaive{1\bar{1}11}(\comring) = \risomLieNaive{11}(\cubel{\comring}{11})$ and $\risomLieNaive{10\bar{1}0}(\conic) = \risomLieNaive{10}(\cubel{\conic}{13})$.
		
		The root gradings of $ L $ induce a \enquote{$5$-grading} $ L = \bigoplus_{i=-2}^2 L_i $ where, for each $ i \in \Set{-2, \ldots, 2} $,
		\[ L_i = \bigoplus_{\alpha \in \Fgrad{i}} L_\alpha = \bigoplus_{\alpha \in \Ggrad{i}} L_\alpha. \]
		Note that
		\begin{gather*}
			L_{-2} = L_{-2,-1} \cong \comring, \qquad L_{2} = L_{2,1} \cong \comring, \qquad L_0 = L_{0,-1} \oplus L_{00} \oplus L_{0,1} \cong \jor \oplus L_{00} \oplus L_{0,-1}, \\
			L_{-1} = L_{-1,-2} \oplus L_{-1,-1} \oplus L_{-1,0} \oplus L_{-1,1} \cong \comring \oplus \jor \oplus \jor \oplus \comring \cong L_{1,-1} \oplus L_{1,0} \oplus L_{1,1} \oplus L_{1,2} = L_1.
		\end{gather*}
		To simplify our notation of elements of $L$, we put $x \defl \risomLieNaive{(-2,-1)}(1)$, $y \defl \risomLieNaive{(2,1)}(1)$, $\adpos{c} \defl \risomLieNaive{(0,1)}(c)$, $\adneg{c} \defl \risomLieNaive{(0,-1)}(c)$ for all $c \in \jor$ and
		\begin{align*}
			[\lambda, c, c', \mu]_- &\defl \risomLieNaive{(-1,-2)}(\lambda) + \risomLieNaive{(-1,-1)}(c) + \risomLieNaive{(-1,0)}(c') + \risomLieNaive{(-1,1)}(\mu), \\
			[\lambda, c, c', \mu]_+ &\defl \risomLieNaive{(1,-1)}(\lambda) + \risomLieNaive{(1,0)}(c) + \risomLieNaive{(1,1)}(c') + \risomLieNaive{(1,2)}(\mu)
		\end{align*}
		for all $\lambda, \mu \in \comring$, $c,c' \in \jor$. (The notation $\adpos{c}$ and $\adneg{c}$ is motivated by the fact that the adjoint map $ \map{\ad}{L}{\End(L)}{\ell}{(r \mapsto [\ell,r])} $ is injective on $L_0$.) We refer to the space $[\comring,J,J,\comring] \defl \Set{[\lambda, c, c', \mu] \given \lambda, \mu \in \comring, c,c' \in J}$ as the \emph{Brown algebra} because, if $\comring$ is a field of characteristic not $2$ or $3$, then it is a special case of a construction given by Brown \cite{Brown}. However, we emphasise that this is mere terminology in our context, and that we do not equip $[\comring,J,J,\comring]$ with any kind of algebra structure.
	\end{miscthm}
	
	\begin{miscthm}[Root homomorphisms]\label{roothom-twisted}
		In practice, it is better to work with the modified root homomorphisms
		\[ \map{\risomLie{\alpha}}{}{}{b}{\risomLieNaive{\alpha}\brackets[\big]{\epsilon_\alpha \gamma_1^{\tau_1(\alpha)} \gamma_2^{\tau_2(\alpha)} \gamma_3^{\tau_3(\alpha)} b}} \]
		for suitably chosen $\epsilon_\alpha \in \compactSet{\pm 1}$, $\tau_1(\alpha), \tau_2(\alpha), \tau_3(\alpha) \in \Set{-1,0,1}$. The reasons for this are explained in \cite[Section~10.3]{DMW26}, and our homomorphisms $(\risomLie{\alpha})_{\alpha \in F_4}$ coincide with the ones in \cite[10.29]{DMW26}.
	\end{miscthm}
	
	\begin{miscthm}[Interface for computations in $L$]
		The package \cjma{} can perform computations in $L$. The elements $ x $, $ y $, $ \xi $, $ \zeta $ are available as constants \texttt{LieX}, \texttt{LieY}, \texttt{LieXi}, \texttt{LieZeta}. The maps $ \adpos{} $, $ \adneg{} $ and $ \dd $ are available as \texttt{adPos}, \texttt{adNeg} and \texttt{dd}.
\begin{verbatim}
gap> cub := CubicEl(a1, 1, 2); t1*LieX; LieXi+2*LieZeta+adPos(cub); dd(cub, t1*cub);
((1)*a1)[12]
(t1)*x
ad^+_{((1)*a1)[12]}+xi+(2)*zeta
dd_{((1)*a1)[12],((t1)*a1)[12]}
\end{verbatim}
		Elements of $ L_1 $ and $ L_{-1} $ are constructed with \texttt{BrownPosEl} and \texttt{BrownNegEl} while \texttt{CubicZero} is $ 0_J $.
\begin{verbatim}
gap> BrownNegEl(Zero(ComRing), CubicZero, cub, t1+t2);
[ 0, 0_J, ((1)*a1)[12], t1+t2 ]_-
\end{verbatim}
		As an alternative way to construct \enquote{basic} elements of $L$, we can use \texttt{LieRootHomF4} to compute the root homomorphisms $(\risomLie{\alpha})_{\alpha \in F_4}$ from~\ref{roothom-twisted}.
\begin{verbatim}
gap> LieRootHomF4([1, 0, 0, 1], a1); LieRootHomF4([2, 0, 0, 0], t1);
[ 0, 0_J, ((g1^-1)*a1)[12], 0 ]_+
(-t1)*y
\end{verbatim}
		The Lie bracket on $L$ is defined in \cite[Table~1]{DMW26} in terms of the structural maps of the cubic norm structure $J$. It can be computed in \cjma{} using the \texttt{*} operator. For example, $[\adpos{c}, \adneg{c'}] = -\dd_{c,c'}$ for all $c,c' \in J$ and $[x,y] = \xi$:
\begin{verbatim}
gap> CubicPosToLieEmb(cub)*CubicNegToLieEmb(cub); LieX*LieY;
(-1)*dd_{((1)*a1)[12],((1)*a1)[12]}
xi
\end{verbatim}
		Just like for elements of $\comring$ and $\conic$, there is a function \texttt{Simplify} for elements of $L$. It applies the aforementioned \texttt{Simplify} functions to the the direct summands $(L_\alpha)_{\alpha \in F_4}$ of $L$ and some additional simplification laws to the summand $L_{0000}$. In the following example, the relations $\contr(a^2)=\contr(a)^2-2\connorm(a)$ for $a \in \conic$ (\ref{conic-identities}~\ref{id:tr-square}) and $\dd_{\cubel{a}{12}, \cubel{b}{23}} = \gamma_2 \dd_{\cubel{1}{11}, \cubel{ab}{13}}$ for $a,b \in \conic$ (\cite[10.10~(iii)]{DMW26}) are used.
\begin{verbatim}
gap> lam := ConicTr(a1^2)-ConicTr(a1)^2+2*ConicNorm(a1);;
gap> l := lam*LieX+dd(CubicEl(a1, 1, 2), CubicEl(a2, 2, 3)); Simplify(l);
(-tr(a1)^2+2*n(a1)+tr(a1a1))*x + dd_{((1)*a1)[12],((1)*a2)[23]}
dd_{(1)[11],((g2)*(a2'*a1'))[31]}
\end{verbatim}
	\end{miscthm}
	
	With the Lie algebra $L$ at hand, we can now study certain automorphisms of $L$.
	
	\begin{miscthm}[Exponential maps on $L$]\label{exp}
		Since $F_4$ is finite, there exists $n \in \IN $ such that $\beta+n\alpha \notin F_4$ for all $\alpha,\beta \in F_4$. An inspection of $F_4$ shows that we can take $n=5$. Hence for any $\alpha \in F_4$ and any $\ell \in L_\alpha$, we have $\ad_{\ell}^5=0$. Thus if $2$ and $3$ are invertible in $\comring$, then we can define an exponential map $\exp(\ell) \defl \sum_{i=0}^4 \frac{1}{i!} \ad_\ell^i$ which is an automorphism of $L$ by standard properties of exponentials. By explicitly computing the action of the exponentials on $L$, we find that all occurrences of $1/2$ and $1/3$ in this action cancel out. For example,
		\[ \brackets[\big]{\exp(x)}(y) = y + [x,y] + \frac{1}{2} \sqbrackets[\big]{x, [x,y]} = y+\xi + \frac{1}{2} [x,\xi] = y+\xi + \frac{1}{2} \cdot (2x) = y+\xi+x  \]
		(where $x,y$ are the fixed elements from~\ref{lie-alg}). Using such explicit formulas, we can define an endomorphism $\exp(\ell)$ of $L$ (as a $\comring$-module) for any $\alpha \in F_4$ and $\ell \in L_\alpha$ without assuming that $2$ and $3$ are invertible (see \cite[Section~7.2]{DMW26}). This lack of additional assumptions on $\comring$ is a key feature of \cite{DMW26}. The drawback of this approach is that we cannot apply any of the standard properties of exponentials: For example, it is not clear that $\exp(\ell)$ is an automorphism of $ L $ (as a Lie algebra) and that $\exp(\ell+\ell') = \exp(\ell) \circ \exp(\ell')$. These properties have to be verified by a computation, which is done in \cite[7.20]{DMW26}
		(in a more general context than cubic Jordan matrix algebras, and hence without using \cjma{}).
		
		Once we have the exponential maps at our disposal, we can define for any $\alpha \in F_4$ a \emph{root group} $\rootgr{\alpha} \defl \exp(L_\alpha)$ and a \emph{root homomorphism} $\map{\risomGrp{\alpha}}{}{}{b}{\exp(\risomLie{\alpha}(b))}$, which is defined on $\comring$ if $\alpha$ is long and on $\conic$ if $\alpha$ is short. The root homomorphisms are available in \cjma{} as \texttt{GrpRootHomF4}, and the concatenation of endomorphisms of $ L $ may be computed with \texttt{*}:
\begin{verbatim}
gap> phi := GrpRootHomF4([-2,0,0,0], One(ComRing));; phi(LieY);
(1)*x + xi + (1)*y
gap> psi := GrpRootHomF4([1,0,0,1], a1);; psi(phi(LieY)) = (psi*phi)(LieY);
true
\end{verbatim}
	(Recall that a double semicolon suppresses the printing of the preceding expression.)
	\end{miscthm}
	
	\begin{miscthm}[Testing equality of automorphisms]\label{test-equal-comm}
		Let $ \phi,\psi $ be two automorphisms of $ L $. We want to decide whether $ \phi=\psi $.
		Recall the definition of the $ 5 $-grading $ F_4 = \bigsqcup_{i=-2}^2 \Fgrad{i} $ from~\ref{G2F4}.
		For any $ a \in \conic $, put
		\[ S(a) \defl \Set{\risomLie{\alpha}(a) \given \alpha \in \Fgrad{-2} \cup \Fgrad{1} \text{ short}} \cup \Set{\risomLie{\alpha}(1_\comring) \given \alpha \in \Fgrad{-2} \cup \Fgrad{1} \text{ long}}. \]
		Since $ L_{-2} \cup L_1 $ generates $ L $ as a Lie algebra, we have $ \phi=\psi $ if and only if for all $ a \in \conic $ and all $ \ell \in S(a) $, we have $ \phi(\ell) = \psi(\ell) $. If $ a $ is \texttt{ConicAlgIndet(i)} for an index $ i $ such that the definition of $ \phi $ and $ \psi $ does not involve \texttt{ConicAlgIndet(i)} (in other words, \enquote{$\phi$, $\psi$ are independent of $a$}), then it suffices to check that $ \phi(\ell) = \psi(\ell) $ for all $ \ell \in S(a) $. For example, if we wanted to prove the equality of \texttt{GrpRootHomF4([1,0,0,1], a1+a2)} and \texttt{GrpRootHomF4([1,0,0,1], a1)*GrpRootHomF4([1,0,0,1], a2)}, then we can do this by evaluation on $S(\mathtt{a3})$ but not by evaluation on $S(\mathtt{a1})$ or $S(\mathtt{a2})$.
		
		Let $ a $ be \texttt{ConicAlgIndet(i)} for the largest possible index $ i $ (as specified initially with \texttt{InitCJMA}) and assume that the definition of $ \phi $ and $ \psi $ does not involve $ a $. The function \texttt{TestEquality} computes $ \phi(\ell)-\psi(\ell) $ for all $ \ell \in S(a) $ and then checks whether all these terms reduce to zero with \texttt{Simplify}. If this is the case, then indeed $ \phi=\psi $ and the function returns \texttt{true}. Otherwise it returns a list that describes the non-zero expressions that were computed in this way, so that the user may check if they can be reduced to zero by hand. For example, consider the root base $(\delta_1, \ldots, \delta_4)$ of $F_4$ given by
		\begin{equation}\label{eq:root-base}
			\delta_1 \defl (1,1,-1,-1), \qquad \delta_2 \defl (-2,0,0,0), \qquad \delta_3 \defl (1,-1,0,0), \qquad \delta_4 \defl (0,1,1,0)
		\end{equation}
		The following computation proves the commutator relation $ [\risomGrp{\delta_2}(t), \risomGrp{\delta_3}(a)] = \risomGrp{\delta_2+\delta_3}(-ta) \risomGrp{\delta_2+2\delta_3}(t \connorm(a)) $:
\begin{verbatim}
gap> d2 := [-2,0,0,0];; d3 := [1,-1,0,0];; d4 := [0,1,1,0];; hom := GrpRootHomF4;;
gap> comm := hom(d2, -t1)*hom(d3, -a1)*hom(d2, t1)*hom(d3, a1);;
gap> TestEquality(comm, hom(d2+d3, -t1*a1) * hom(d2+2*d3, t1*ConicNorm(a1)));
true
\end{verbatim}
		We now prove $ [\risomGrp{\delta_4}(a), \risomGrp{\delta_3}(b)] = \risomGrp{\delta_3+\delta_4}(\conconj{b} \conconj{a}) $. The following computation tells us that the two sides of this equation differ by $ \dd(\cubel{1}{22}, \cubel{-\gamma_3^{-1} \connorm(b) a + \gamma_3^{-1} (ab) \conconj{b}}{12}) $ when evaluated on $ x $.
\begin{verbatim}
gap> comm := hom(d4, -a1)*hom(d3, -a2)*hom(d4, a1)*hom(d3, a2);;
gap> TestEquality(comm, hom(d3+d4, ConicInv(a2)*ConicInv(a1)));
[ [ (1)*x, dd_{(1)[22],((-n(a2)/g3)*a1+(1/g3)*((a1*a2)*a2'))[12]} ] ]
\end{verbatim}
		Since $ \connorm(b)a = a \connorm(b) = a \connorm(\conconj{b}) = (ab) \conconj{b} $ by~\ref{conic-identities}, we conclude that the desired commutator relation holds.
	\end{miscthm}
	
	Finally, we connect the root groups $(\rootgr{\alpha})_{\alpha \in F_4}$ to the theory of root graded groups.
	
	\begin{miscthm}[Root graded groups]\label{rgg}
		Let $\roots$ be a finite root system. A \emph{$\roots$-graded group} is a group $G$ together with a family $(\rootgrvar{\alpha})_{\alpha \in \roots}$ of subgroups of $G$ satisfying certain axioms, most importantly:
		\begin{enumerate}[(i)]
			\item $[\rootgrvar{\alpha}, \rootgrvar{\beta}] \subseteq \gen{\rootgrvar{\gamma} \given \gamma \in \roots, \gamma = i\alpha+j\beta \text{ for some } i,j \in \IN_{\ge 1}}$ for all non-proportional $ \alpha, \beta \in \roots $.
			
			\item \label{rgg:weyl}For each $ \alpha \in \roots $, there exists an \emph{$\alpha$-Weyl element}: An element $ w_\alpha \in \rootgrvar{-\alpha} \rootgrvar{\alpha} \rootgrvar{-\alpha} $ such that $ \rootgrvar{\beta}^{w_\alpha} = \rootgrvar{\beta^{\reflbr{\alpha}}} $ for all $ \beta\in \roots $.
		\end{enumerate}
		For a precise definition of root graded groups, which generalise Jacques Tits' notion of RGD-systems, see \cite[2.5.2]{WiedemannPhD}. One of the main results of \cite{WiedemannPhD} is that for any $F_4$-graded group $(G, (\rootgrvar{\alpha})_{\alpha \in F_4})$, there exist a commutative ring $\comring'$ and a conic $\comring'$-algebra $\conic'$ which \emph{coordinatise} $G$ (see \cite[10.7.8]{WiedemannPhD}). This means that there exist an isomorphism $\map{\risomGrp{\alpha}}{\comring'}{\rootgrvar{\alpha}}{}{}$ for each long $\alpha \in F_4$ and an isomorphism $\map{\risomGrp{\beta}}{\conic'}{\rootgrvar{\beta}}{}{}$ for each short $\beta \in F_4$ such that certain explicit commutator relations involving the structural maps of $\conic$ are satisfied. The commutator relations that we have verified in~\ref{test-equal-comm} are examples of the required relations.
	\end{miscthm}
	
	\begin{miscthm}[The existence problem for $F_4$-graded groups]\label{ex-prob}
		The coordinatisation result for $F_4$-graded groups in~\ref{rgg} gives rise to the following existence question: Given $\comring$ and $\conic$, does there exist an $F_4$-graded group that is coordinatised by $\conic$? This question was the initial motivation for the work in \cite{DMW26}. It is shown in \cite[11.15]{DMW26} that the group $G$ generated by the root groups $(\rootgr{\alpha})_{\alpha \in F_4}$ from \ref{exp} is indeed $ F_4 $-graded and coordinatised by $ \conic $, so the question has a positive answer. We have already seen in~\ref{test-equal-comm} how to verify some commutator relations in $G$, and the remaining commutator relations can be verified with \cjma{} in a similar way.
		
		The biggest computational effort, however, lies in the verification of Axiom~\ref{rgg}~\ref{rgg:weyl}. Let $\rootbase$ be the root base of $F_4$ from \eqref{eq:root-base}. By \cite[2.2.7]{WiedemannPhD}, it suffices to show that there exists a $\delta$-Weyl element for each $\delta \in \rootbase$. Let $\delta \in \rootbase$ be arbitrary and put $w_\delta \defl \risomGrp{-\delta}(-1) \risomGrp{\delta}(1) \risomGrp{-\delta}(-1) \in G$. We show that $w_\delta$ is a $\delta$-Weyl element by verifying with \cjma{} (and a few final steps by hand) that
		\begin{align*}
			\risomGrp{\alpha}(a)^{w_\delta} \in \Set{\risomGrp{\alpha^{\reflbr{\delta}}}(a), \risomGrp{\alpha^{\reflbr{\delta}}}(-a), \risomGrp{\alpha^{\reflbr{\delta}}}(\conconj{a}), \risomGrp{\alpha^{\reflbr{\delta}}}(-\conconj{a})} \quad \text{and} \quad \risomGrp{\beta}(t)^{w_\delta} \in \Set{\risomGrp{\beta^{\reflbr{\delta}}}(t), \risomGrp{\beta^{\reflbr{\delta}}}(-t)}
		\end{align*}
		for all short $\alpha \in F_4$, long $\beta \in F_4$ and all $t \in \comring$, $a \in \conic$. This requires up to four equality tests of automorphisms (as in~\ref{test-equal-comm}) for each short $ \alpha \in F_4 $ and up to two such tests for each long $ \beta \in F_4 $. In total, we have to perform $ \abs{F_4} \cdot \abs{\rootbase} = 48 \cdot 4 $ such computations to show that $ \rootgr{\gamma}^{w_\delta} = \rootgr{\gamma^{\reflbr{\delta}}} $ for all $ \gamma \in F_4 $, $ \delta \in \rootbase $.
	\end{miscthm}
	
	\section{Implementation details}\label{sec:implement}
	
	\begin{miscthm}[The underlying rings]\label{imp:rings}
		Denote by $ m_1 $, $ m_2 $, $ m_3 $ the first three arguments passed to \texttt{InitCJMA} during initialisation. Internally, \cjma{} sets up a polynomial ring \texttt{ComRing} with $ m_1+m_2+3+\abs{T(m_2,m_3)} $ indeterminates (where the set $ T(m_2,m_3) $ will be defined momentarily) and a free nonassociative \texttt{ComRing}-algebra \texttt{ConicAlg} on $ 2m_2 $ generators. These rings represent $ \comring $ and $ \conic $, respectively, and their generators are named in the following manner. The $ 2m_2 $ generators of \texttt{ConicAlg} are named \texttt{a1}, \dots, \texttt{a}$ m_2 $, \texttt{a1'}, \dots, \texttt{a}$ m_2 $\texttt{'}. The first $ 3+m_1+m_2 $ indeterminates in \texttt{ComRing} are the symbols \texttt{g1}, \texttt{g2}, \texttt{g3}, \texttt{t1}, \dots, \texttt{t}$ m_1 $, \texttt{n(a1)}, \ldots, \texttt{n(a}$ m_2 $\texttt{)}. The remaining $ \abs{T(m_2,m_3)} $ indeterminates in \texttt{ComRing} are described as follows.

		For any set $ S $ of symbols, denote by $ \gen{S} $ the set of nonassociative words over $ S $ (that is, the free magma generated by $ S $) and by $ \gen{S}_{\le \ell} $ the set of nonassociative words of length at most $ \ell \in \IN $. Let $ A_{m_2} $ and $ A_{m_2}' $ denote the sets of symbols \texttt{a1}, \dots, \texttt{a}$ m_2 $ and \texttt{a1'}, \dots, \texttt{a}$ m_2 $\texttt{'}, respectively. For each $ d \in \gen{A_{m_2} \cup A'_{m_2}}_{\le m_3} $, we need an indeterminate in \texttt{ComRing} to represent $ \contr(d) $. However, many of these trace values actually coincide:
		Indeed, we deduce from~\ref{conic-identities}~\ref{id:tr} and~\ref{conic-identities}~\ref{id:conj} that
		\begin{align}\label{eq:tr-indet}
			\contr(abc) = \contr(bca) = \contr(cab) = \contr(\conconj{c} \conconj{b} \conconj{a}) = \contr(\conconj{b} \conconj{a} \conconj{c}) = \contr(\conconj{a} \conconj{c} \conconj{b}) \quad \text{and} \quad \contr\brackets[\big]{a \conconj{a}b} = \contr\brackets[\big]{\connorm(a)b} = \connorm(a) \contr(b)
		\end{align}
		for all $ a,b,c \in \conic $. Hence we can find a proper subset $ T(m_2,m_3) $ of $ \gen{A_{m_2} \cup A'_{m_2}}_{\le m_3} $ such that for each $ \tilde{d} \in \gen{A_{m_2} \cup A'_{m_2}}_{\le m_3} $, there exist $ d \in T(m_2,m_3) $, $ e \in \gen{A_{m_2}}_{\le m_3} $ with $ \contr(\tilde{d}) = \connorm(e)\contr(d) $ (where $\connorm(e)$ should be read as $1_\comring$ if $e$ is the empty word). We introduce an indeterminate \texttt{tr(d)} in $ \comring $ for each $ d \in T(m_2,m_3) $.
		
		Alternatively, it would be possible to deal with the relations~\eqref{eq:tr-indet} only later in the \texttt{Simplify} procedure, but this would greatly increase the number of trace indeterminates that have to be set up in \texttt{ComRing}.
	\end{miscthm}
	
	\begin{miscthm}[The remaining structures]
		The cubic Jordan matrix algebra $J$, the Brown algebra $[\comring, J, J, \comring]$, the Lie algebra $L$ and its automorphism group $\Aut(L)$ are constructed using \texttt{ArithmeticElementCreator}, a tool from the GAP library which automatically defines the operators \texttt{+} and \texttt{*} for such objects. The internal representation of these elements is straightforward in most cases: For examples, elements of $[\comring, J, J, \comring]$ may be represented as $4$-tuples. An element $t_1 \dd_{c_1,d_1}+ t_2 \dd_{c_2,d_2} + \cdots$ of $L_{00}$ (with $t_i \in \comring$, $c_i,d_i \in J$) is represented as a list \texttt{[[t1, c1, d1], [t2, c2, d2], ...]}. Of course, such a representation is not unique---this problem is alleviated to some extent by the \texttt{DDSanitizeRep} function, which applies basic identities such as $t_1\dd_{c,d}+t_2\dd_{c,d}=(t_1+t_2)\dd_{c,d}$, and by the more sophisticated \texttt{Simplify} procedure. The additional elements $ \xi $ and $ \zeta $ of $ L_{00} $ are handled in a similar way. Elements of $\Aut(L)$ are represented as GAP functions from $L$ to $L$.
	\end{miscthm}
	
	\begin{miscthm}[Simplification on the underlying rings]
		The \texttt{Simplify} routine for \texttt{ConicAlg} relies on~\ref{conic-identities}~\ref{id:conj} to replace occurences of $a+\conconj{a}$ by $\contr(a)1_\conic$. This substitution is performed by the subroutine \texttt{MakeTraces}.
		
		We keep the notation of~\ref{imp:rings}. The \texttt{Simplify} routine on \texttt{ComRing} performs the following substitutions:
		\begin{enumerate}
			\item Using the definition of the conjugation on $ \conic $, we see the following: For any $ \ell \in \IN $ and $b \in \gen{A_{m_2} \cup A'_{m_2}}_{\le \ell}$, we can substitute $\contr(b)$ by a sum of products of terms in $ \Set{\contr(c) \given c \in \gen{A_{m_2}}_{\le \ell}} $. For example,
			\[ \contr(a_i \conconj{a_j}) = \contr\brackets[\big]{a_i \brackets{\contr(a_j) 1_\conic -a_j}} = \contr(a_i) \contr(a_j) - \contr(a_i a_j) \]
			for all $ i \ne j \in \Set{1,\ldots,m_2} $, so we may substitute $ \contr(a_i \conconj{a_j}) $ by $ \contr(a_i) \contr(a_j) - \contr(a_i a_j) $.
			The \texttt{Simplify} function on \texttt{ComRing} performs such substitutions for $ \ell \in \Set{2,3} $ using GAP's \texttt{Value} function for substitution in polynomials. We could also perform similar substitutions for higher values of $ \ell $, but this is not necessary for our applications.
					
			\item Each occurrence of $ \contr(a_i a_i) $ for $ i \in \Set{1,\ldots,m_2} $ is substituted by $ \contr(a_i)^2 - 2\connorm(a_i) $ (see~\ref{conic-identities}~\ref{id:tr-square}).
		\end{enumerate}
	\end{miscthm}
	
	\begin{miscthm}[Simplification on the remaining structures]
		The \texttt{Simplify} routines on $J$ and $[k,J,J,k]$ simply call those of \texttt{ComRing} and \texttt{ConicAlg} on the components of the structure in question. The \texttt{Simplify} routine on $L$ does the same, except that the direct summand $L_{00}$ needs a more sophisticated approach. We use the following relations to simplify the representation of an element of $L_{00}$:
		\begin{enumerate}
			\item Let $i,j \in \Set{1,2,3}$ and put $Z_{i \to j} \defl Z_{i1,1j} + Z_{i2,2j}  + Z_{i3,3j}$ where $Z_{ij,pq} \defl \gen{\dd_{c,c'} \given c \in J_{ij}, c' \in J_{pq}}_\comring$ for $i,j,p,q \in \Set{1,2,3}$. For $i \ne j$, it is shown in \cite[10.14]{DMW26} that $\map{\phi_{i \to j}}{\conic}{Z_{i \to j}}{a}{\dd_{\cubel{1}{ii},\cubel{a}{ij}}}$ is an isomorphism of $\comring$-modules. The key ingredient in the proof is the relation \cite[10.10~(iii)]{DMW26}:
			\begin{align*}
				\dd(\cubel{a}{jl},\cubel{b}{li}) = \dd(\cubel{1_\conic}{jl}, \cubel{ab}{li}) = \gamma_l \dd(\cubel{1}{jj}, \cubel{ab}{ji}) = \gamma_l \dd(\cubel{ab}{ji}, \cubel{1}{ii})
			\end{align*}
			We can use this relation to rewrite an element of $Z_{i \to j}$ as $\phi_{i \to j}(a)$ for a unique $a \in \conic$.
			
			\item By the $\comring$-bilinearity of $\map{}{}{}{(c,c')}{\dd_{c,c'}}$, any element of $Z_{ii,ii}$ for $i \in \Set{1,2,3}$ can be rewritten as $\dd(\cubel{1}{ii},\cubel{t}{ii})$ for a unique $t \in \comring$.
			
			\item Using the Peirce decomposition from \cref{ex:peirce}, many trivial terms can be removed, such as $\dd_{\cubel{a}{12}, \cubel{t}{33}}$. Of course, this substitution requires us to first decompose any summand $\dd_{c,c'}$ for $c,c' \in J$ as a sum of terms $\dd_{\cubel{a}{ij},\cubel{b}{pq}}$ for $a,b \in \conic \cup \comring$, $i,j,p,q \in \Set{1,2,3}$.
			
			\item 
			By \cite[10.10]{DMW26}, we can replace $\dd_{\cubel{1}{33}, \cubel{1}{33}}$ by $2\zeta - \xi - \dd_{\cubel{1}{11}, \cubel{1}{11}} - \dd_{\cubel{1}{22}, \cubel{1}{2}}$ .
			
			\item By \cite[10.10~(iv)]{DMW26}, we can replace $\dd_{\cubel{a}{ij},\cubel{\conconj{a}}{ji}}$ for $a \in \conic$ and $i \ne j \in \Set{1,2,3}$ by
			\[ \gamma_i \gamma_j \connorm(a) \brackets[\big]{\dd_{\cubel{1}{ii}, \cubel{1}{ii}} + \dd_{\cubel{1}{jj}, \cubel{1}{jj}}}. \]
			We also use the linearised version of this identity,
			\[ \dd_{\cubel{a}{ij}, \cubel{\conconj{b}}{ji}} + \dd_{\cubel{b}{ij}, \cubel{\conconj{a}}{ji}} = \gamma_i \gamma_j \contr(a \conconj{b}) \brackets[\big]{\dd_{\cubel{1}{ii}, \cubel{1}{ii}} + \dd_{\cubel{1}{jj}, \cubel{1}{jj}}}, \]
			as follows: Any sum $\lambda_1\dd_{\cubel{a}{ij}, \cubel{\conconj{b}}{ji}} + \lambda_2\dd_{\cubel{b}{ij}, \cubel{\conconj{a}}{ji}}$ for $\lambda_1, \lambda_2 \in \comring$ is replaced by 
			\[ (\lambda_1 - \lambda_2) \dd_{\cubel{a}{ij}, \cubel{\conconj{b}}{ji}} + \lambda_2 \brackets[\big]{\dd_{\cubel{1}{ii}, \cubel{1}{ii}} + \dd_{\cubel{1}{jj}, \cubel{1}{jj}}}. \]
		\end{enumerate}
	\end{miscthm}

	\section{Comparison to the literature}\label{sec:literature}
	
	We are not aware of any literature or software that discusses symbolic computation in conic algebras.
	However, there has been some work on the related problem of symbolic computation in alternative algebras.
	Recall from \cref{sec:intro} the problem of deciding whether an element of a free nonassociative $\comring$-algebra is an identity in alternative $\comring$-algebras.
	The computer algebra system \emph{Albert} \cite{Albert-Algo,Albert-System,Albert-Progress,Albert-Homepage} solves this problem (and similar problems for arbitrary classes of nonassociative algebras defined by certain identities) under the restriction that the ground ring $\comring$ is a field of characteristic not $2$. The assumption that $2$ is invertible in $\comring$ is necessary because Albert works with the linearised alternative laws
	\begin{equation}\label{}
		(aa')b + (a'a)b = a(a'b) + a'(ab) \quad \text{and} \quad (ab)b' + (ab')b = a(bb') + a(b'b) \quad \text{for all } a,a',b,b' \in \alt,
	\end{equation}
	which are equivalent to the usual alternative laws~\eqref{eq:alternative} under this assumption. The assumption that $\comring$ is a field is necessary because Albert works by solving linear equation systems. Further, for implementation reasons, Albert only supports prime fields $\comring$ of order at most $251$.
	
	A different approach to this decision problem, which uses term rewriting, is taken in \cite{Widiger-Alternative}. It works over $\comring=\IZ$, but is only guaranteed to be complete if all \enquote{monomials} in the identity $f$ that we want to verify have total degree at most $4$. Thus if $f$ is an identity that contains a higher-degree monomial, then this approach might fail to detect that $f$ is an identity. However, we are not aware of any publicly available software implementation of the ideas in \cite{Widiger-Alternative}.
	
	\section*{Acknowledgements}
	
	I want to thank Till Eisenbrand, Max Horn and Martin Wagner for several helpful discussions during our working group seminar which improved the presentation of this paper. This work was partially supported by a DAAD postdoctoral research scholarship and by the Deutsche Forschungsgemeinschaft (DFG, German Research Foundation), Project-ID 286237555, TRR 195.
	
	\bibliographystyle{alpha}
	\bibliography{CubicJordanMatrixAlg-bib}

@article {Albert-Algo,
    AUTHOR = {Hentzel, Irvin Roy and Jacobs, David Pokrass},
     TITLE = {A dynamic programming method for building free algebras},
   JOURNAL = {Comput. Math. Appl.},
  FJOURNAL = {Computers \& Mathematics with Applications. An International
              Journal},
    VOLUME = {22},
      YEAR = {1991},
    NUMBER = {12},
     PAGES = {61--66},
      ISSN = {0898-1221,1873-7668},
       DOI = {10.1016/0898-1221(91)90148-W},
}

@incollection {Albert-System,
    AUTHOR = {Jacobs, Pokrass and Muddana, Sekhar V. and Offutt, A.
              Jefferson},
     TITLE = {A computer algebra system for nonassociative identities},
 BOOKTITLE = {Hadronic mechanics and nonpotential interactions, {P}art 1
              ({C}edar {F}alls, {IA}, 1990)},
     PAGES = {185--195},
 PUBLISHER = {Nova Sci. Publ., Commack, NY},
      YEAR = {1992},
      ISBN = {1-56072-035-2},
}

@inproceedings{Albert-Progress,
author = {Jacobs, David P.},
title = {The {A}lbert nonassociative algebra system: a progress report},
year = {1994},
isbn = {0897916387},
publisher = {Association for Computing Machinery},
address = {New York, NY, USA},
url = {https://doi.org/10.1145/190347.190358},
doi = {10.1145/190347.190358},
booktitle = {Proceedings of the International Symposium on Symbolic and Algebraic Computation},
pages = {41--44},
numpages = {4},
location = {Oxford, United Kingdom},
series = {ISSAC '94}
}

@misc{Albert-Homepage,
	title = {Albert's page},
	author = {Jacobs, David P.},
	howpublished= {\url{https://dpj.people.clemson.edu/albertstuff/albert.html}},
	note = {Access date: 2026-04-01}
}

@manual{GAP4,
organization = "The GAP~Group",
title        = "{GAP -- Groups, Algorithms, and Programming,
                Version 4.15.1}",
year         = 2025,
url          = "\url{https://www.gap-system.org}",
}

@misc{DMW26,
      title={From cubic norm pairs to {$G_2$}- and {$F_4$}-graded groups and {L}ie algebras}, 
      author={De Medts, Tom and Wiedemann, Torben},
      year={2026},
      eprint={2602.06147},
      archivePrefix={arXiv},
      primaryClass={math.RA},
      note={Preprint, 114 pages, \url{https://arxiv.org/abs/2602.06147}}, 
}

@book{WiedemannPhD,
	author		={Wiedemann, Torben},
	title		={Root Graded Groups},
	year		={2024},
	note		={PhD Thesis, \url{https://doi.org/10.22029/jlupub-18373}},
	publisher	={Justus Liebig University Giessen},
}

@misc{CJMA,
	title={{C}ubic{J}ordan{M}atrix{A}lg},
	author={Wiedemann, Torben},
	year={2026},
	note={GAP package, version 1.0.4, \url{https://doi.org/10.5281/zenodo.19221614}}
}

@article {Brown,
    AUTHOR = {Brown, Robert B.},
     TITLE = {A new type of nonassociative algebras},
   JOURNAL = {Proc. Nat. Acad. Sci. U.S.A.},
  FJOURNAL = {Proceedings of the National Academy of Sciences of the United
              States of America},
    VOLUME = {50},
      YEAR = {1963},
     PAGES = {947--949},
      ISSN = {0027-8424},
       DOI = {10.1073/pnas.50.5.947},
}

@book {GPR24,
    AUTHOR = {Garibaldi, Skip and Petersson, Holger P. and Racine, Michel
              L.},
     TITLE = {Albert algebras over commutative rings---the last frontier of
              {J}ordan systems},
    SERIES = {New Mathematical Monographs},
    VOLUME = {48},
 PUBLISHER = {Cambridge University Press, Cambridge},
      YEAR = {2024},
     PAGES = {xxii+655},
      ISBN = {978-1-009-42685-5},
}

@article {McCrimmon66,
    AUTHOR = {McCrimmon, Kevin},
     TITLE = {A general theory of {J}ordan rings},
   JOURNAL = {Proc. Nat. Acad. Sci. U.S.A.},
  FJOURNAL = {Proceedings of the National Academy of Sciences of the United
              States of America},
    VOLUME = {56},
      YEAR = {1966},
     PAGES = {1072--1079},
      ISSN = {0027-8424},
       DOI = {10.1073/pnas.56.4.1072},
}

@book {Schafer,
    AUTHOR = {Schafer, Richard D.},
     TITLE = {An introduction to nonassociative algebras},
      NOTE = {Corrected reprint of the 1966 original},
 PUBLISHER = {Dover Publications, Inc., New York},
      YEAR = {1995},
     PAGES = {x+166},
      ISBN = {0-486-68813-5},
   MRCLASS = {17-01},
  MRNUMBER = {1375235},
}

@incollection {Widiger-Alternative,
    AUTHOR = {Widiger, Alfred},
     TITLE = {Deciding degree-four-identities for alternative rings by
              rewriting},
 BOOKTITLE = {Symbolic rewriting techniques ({A}scona, 1995)},
    SERIES = {Progr. Comput. Sci. Appl. Logic},
    VOLUME = {15},
     PAGES = {277--288},
 PUBLISHER = {Birkh\"auser, Basel},
      YEAR = {1998},
      ISBN = {3-7643-5901-3},
   MRCLASS = {68Q42 (17D05 68T15)},
  MRNUMBER = {1624604},
MRREVIEWER = {Ralph\ W.\ Wilkerson},
       DOI = {10.1007/978-3-0348-8800-4\_14},
       URL = {https://doi.org/10.1007/978-3-0348-8800-4_14},
}
\end{document}